\documentclass{elsarticle}
\usepackage{lineno,hyperref}
\usepackage{upgreek}
\modulolinenumbers[5]
\usepackage{array}
\usepackage{amssymb,amsmath,amsthm}
\usepackage{epsfig,multicol,graphicx,epstopdf}
\usepackage{color}
\usepackage{geometry}
\journal{\textcolor{red}{Journal Name}}
\usepackage{graphicx}
\hypersetup{colorlinks,linkcolor={blue},citecolor={blue},urlcolor={red}} 
\newtheorem{The}{Theorem}

\begin{document}

\begin{frontmatter}

\title{Lie group theory for nonlinear fractional K(m,n) type equation with variable coefficients}

\author[a]{H. Jafari}%\corref{corA}
\ead{jafarh@unisa.ac.za}
\author[b]{N. Kadkhoda}
\ead{kadkhoda@buqaen.ac.ir}

\author[c,d]{D. Baleanu \corref{corA}} 
\cortext[corA]{Corresponding author.}
\ead{dumitru@cankaya.edu.tr}
\address[a] {Department of Mathematical Sciences, University of South Africa, UNISA0003, South Africa}
\address[b]{Department of Mathematics, Faculty of Basic Sciences, Bozorgmehr University Of Qaenat,Qaenat,   Iran}
\address[c] {Department of Mathematics, Faculty of Art and Sciences,Cankaya University, Ankara, Turkey}
\address[d]{Institute of Space Sciences, Magurele-Bucharest,Romania}

\begin{abstract}
We investigated the analytical solution of fractional order K(m,n) type equation with variable coefficient which is  an extended type of KdV equations into a genuinely nonlinear dispersion regime.\\
By using the Lie symmetry analysis, we obtain the Lie point symmetries for this type of  time-fractional partial differential equations (PDE). Also 
we present  the corresponding  reduced   fractional  differential equations (FDEs) corresponding to the time-fractional K(m,n) type equation.
\end{abstract}

\end{frontmatter}

% ----------------------------------
\section{Introduction}\label{sec:1}
Most problems in engineering, biology, applied mathematics and physics  might be better modeled  by using ordinary/partial differential equations with fractional (arbitrary) order. The method of group analysis for ordinary/partial differential equations, originally advocated by the Norwegian mathematician Sophus Lie during 1870s.
The tangent structural equations under transformation groups is the  fundamental  idea of symmetry analysis. 
Numerous methods developed to solve differential equations based on Lie symmetry analysis.\\
In last few decades, many researcher studied different class of the fractional partial differential equations (FPDE).
These equations arise in Various branches of sciences such as physics, biology, viscoelastic materials, electrochemistry, signal processing, fluid mechanics
\cite{Kilb,Pod, Bau206,Bau100,Bau2300}. Integrals and derivatives are of any order in the fractional calculus \cite{Kilb,Pod}.  In the recent years, finding exact solutions of FDEs  has gained much attention.

Many researchers have presented various techniques and methods for obtaining the numerical and analytical solutions of FDE, such as variational
iteration method \cite{He}, the separating variables method \cite{Chen}, operational matrices \cite{jaff},
the fractional complex transform \cite{He1}, the first integral method \cite{Lu},  and so on.
In many years ago, there are many articles to obtain  the analytical solution/s of nonlinear PDE using Lie group theory.
It is important to know, however,   that few of them involve  FDEs\cite{ckwa,Gazizov, jaff20,Kasatk, Liu,  Wei, Wu1}.
Furthermore, already few articles done in symmetries of variable coefficients FDEs  such as \cite{Gaur, Luka}.
Our purpose is to study the time-fractional K(m,n) equation:
\begin{equation}\label{o1}
\frac{\partial^\alpha u}{\partial
t^\alpha}+\zeta (u^m)_x+g(t)(u^n)_{xxx}=0, \quad t>0, \, 0<\alpha\leq 1,
\end{equation}
or  equivalently
\begin{equation}\label{}
D^\alpha_tu+\zeta (u^m)_x+g(t)(u^n)_{xxx}=0,\quad t>0, \, 0<\alpha\leq 1.
\end{equation}
here $m$ and $n\neq0$ are arbitrary constants, $\zeta=\pm 1$ and $g(t)$ is  an arbitrary nonvanishing function of the variable $t$. This equation for $\alpha=1$ and also with constant coefficient for $0<\alpha< 1$ has been discussed in \cite{Chara} and \cite{ne1}.\\ In the follow, we study the above equation with
$m=2, n=3$. Eq. \eqref{o1} is called the K(m,n) equation, when $\zeta=g(t)=1$. Rosenau introduced this equation in 1998 \cite{Rose,Z1} which is described the process of interpretation the role
nonlinear dispersion in the formation of structures in liquid drops.
 \\
This article is organized as follows. In the next Section, it is given the analysis of Lie Symmetry group for a FPDE.  Then in section 3, using Lie group, the Lie point symmetries of equation \eqref{o1} are obtained. In Section 4,  we perform Lie group on the equation \eqref{o1} for obtaining  invariant solutions and reduced fractional ODEs .  Conclusions are given in the Section 5.

\section{ Lie symmetry analysis method for  FPDEs}
%%%%%%%%%%%%%%%%%%%%%%%%%%%%%%%%%%%%%
 According to  the importance of FPDEs  in mathematics and physics,  finding the exact solutions for these
 equations is very important. Although nonlinear FPDEs are difficult to
solve, but many papers have been  presented by scientists.
Studying   differential equations using  the   fundamental method of the  Lie
symmetries  is interesting for many researchers.
In the past century, many researchers have studied in the field of
the Lie groups. Some of them are Ovsiannikov \cite{Ovsy}, Olver
\cite{olv}, Ibragimov \cite{Ibrag}, Baumann \cite{Bauman} and Bluman
\cite{Bluma}.
In this section, finding infinitesimal functions of FPDEs are given.
Let us consider the below form of FPDEs::
\begin{equation}\label{aa}
\mathcal{D}^\alpha_t u=F(x,t,u,u_{(1)},\ldots), \quad \alpha>0.
\end{equation}
Where $D^\alpha_t$   fractional derivative in the sense of Riemann-Liouville \cite{jaff20} and $u$ is depend to  $x,t$.

Similar discussion of PDEs\cite{bk,olv}, we can write
\begin{equation}\label{ac}
\mathcal{D}^\alpha_{\bar{t}}\bar{u}=\mathcal{D}^\alpha_tu+\varepsilon[\eta^{(\alpha)}_{\,\,t}(t,x,u,u_{(1)},u_{(\alpha)},\ldots)]+o(\varepsilon^2).
\end{equation}
In view of by the prolongation formula, for $\eta^{(\alpha)}_{\,\,t}$ we have
\cite{Gazizov}
\begin{equation}\label{ad}
\eta^{(\alpha)}_{\,\,t}=\mathcal{D}^\alpha_t(\eta)+\xi^x \mathcal{D}^\alpha_t(u_x)-\mathcal{D}^\alpha_t(\xi^x
u_x)+\mathcal{D}^\alpha_t(D_t(\xi^t)u)-\mathcal{D}^{\alpha+1}_t(\xi^t u)+\xi^t
D^{\alpha+1}_t u,
\end{equation}
and the total derivative operator $D_t$ is defined by
\begin{equation}
D_t=\frac{\partial}{\partial t}+u_t \frac{\partial}{\partial
u}+u_{xt} \frac{\partial}{\partial u_x}+u_{tt}
\frac{\partial}{\partial u_t}+u_{xxt} \frac{\partial}{\partial
u_{xx}}+\cdots.
\end{equation}
Simplifying \eqref{ad} using the Leibniz formula \cite{ang}
\begin{equation}\label{ae}
\mathcal{D}^\alpha_t[f(t) g(t)]=\sum^\infty_{n=0}\frac{(-1)^{n-1}\alpha
\Gamma(n-\alpha)}{\Gamma(1-\alpha)\Gamma(n+1)}\mathcal{D}^{\alpha-n}_t f(t)
D^n_t g(t),\quad \alpha>0,
\end{equation}
we can write \cite{jeff}:
\begin{eqnarray}\label{af}
\eta^{(\alpha)}_{\,\,t}&=&\frac{\partial^\alpha \eta}{\partial t^\alpha}+(\eta_u-\alpha D_t(\xi^t))\frac{\partial^\alpha
u}{\partial t^\alpha}-u \frac{\partial^\alpha \eta_u}{\partial t^\alpha}+\sum^\infty_{m=1}[\binom{\alpha}{m}\frac{\partial^m
(\eta_u)}{\partial
t^m}\nonumber\\&-&\binom{\alpha}{m+1}D^{m+1}_t(\xi^t)]\mathcal{D}^{\alpha-m}_t(u)
-\sum^\infty_{m=1}\binom{\alpha}{m}\mathcal{D}^{\alpha-m}_t(u_x)D^m_t(\xi^x).
\end{eqnarray}
To obtain coefficients of $X$, we must have:
\begin{equation}\label{ag}
X^{(\alpha)}[\mathcal{D}^\alpha_t u-F(t,x,u,u_{(1)},\ldots)]_{_{\mathcal{D}^\alpha_t
u=F(t,x,u_{(1)},\ldots)}}=0,
\end{equation}
where
\begin{eqnarray}\label{ah}
X^{(\alpha)}&=&\xi^x(t,x,u)\frac{\partial}{\partial
x}+\xi^t(t,x,u)\frac{\partial}{\partial
t}+\eta(t,x,u)\frac{\partial}{\partial u}
+\eta^{(1)}_i(t,x,u,u_{(1)})\frac{\partial}{\partial u_i}+\cdots
\nonumber
\\&+&\eta^{(k)}_{i_1i_2\cdots, i_k}(t,x,u,u_{(1)},\ldots
u_{(k)})\frac{\partial}{\partial u_{i_1i_2\cdots, i_k}}+
\eta^{(\alpha)}_t(t,x,u,\ldots,u_{(\alpha),\cdots})\frac{\partial}{\partial
u^{(\alpha)}_t}.
\end{eqnarray}
Using these relations, we obtain the Lie symmetries.
%%%%%%%%%%%%%%%%%%%%%%%%%%%%%%%%%%
\section{Fractional Lie symmetries  for time-fractional K(2,3)} 
%%%%%%%%%%%%%%%%%%%%%%%%%%%%%%%%%%%%
Now we obtain the infinitesimal generator of the
time-fractional K(2,3) equation
\begin{equation}\label{ai}
\frac{\partial^\alpha u}{\partial
t^\alpha}+(u^2)_x+g(t)(u^3)_{xxx}=0, \quad t>0, \,\, 0<\alpha<1.
\end{equation}
\begin{The}Lie symmetries for Eq. \eqref{ai}, which those are solutions  of determining
equations depend on the selection of the function $g(t)$,  are
\begin{description}
  \item[Case 1:] $0<\alpha<1,\,\alpha\neq\frac{1}{2}\,,\frac{1}{3},\ k,b\neq0$.
\begin{description}
  \item[Case 1.1:]  $g(t)$ be an non-vanishing arbitrary function.\\
   In this case, the infinitesimal generator is given by
\begin{equation}
X_{1.1}=\frac{\partial}{\partial x}.
\end{equation}
For  mentioned g(t) as follows, we have additional symmetries.
  \item[Case 1.2:] $ g(t)=\,kt^b$.\\
  For this case, we have
\begin{equation}
X_{1.2.1}=\frac{\partial}{\partial x}, \,\,
X_{1.2.2}=-t\frac{\partial}{\partial t}+(\alpha-b)x\frac{\partial}{\partial x}+(2\alpha-b)u\frac{\partial}{\partial
u}.
\end{equation}

  \item[Case 1.3:] $ g(t)=\,k$.\\
  In this case, the infinitesimal generators are as follows
\begin{equation}
X_{1.3.1}=\frac{\partial}{\partial x}, \,\, X_{1.3.2}=\alpha
x\frac{\partial}{\partial x}-t\frac{\partial}{\partial t}+2\alpha
u\frac{\partial}{\partial u}.
\end{equation}

\end{description}

 \item[Case 2:] $\alpha=\frac{1}{2},\,\ k,b\neq0$.\\
For $\alpha=\frac{1}{2}$, functions of   $g(t)$  can be obtained as
follows
\begin{equation*}
g(t)=ke^{b t},\,\,kt^b,\, \, k.
\end{equation*}
\begin{description}

\item[Case 2.1:]  $g(t)= ke^{bt}$.\\
   In this case, the infinitesimal generator is given by
\begin{equation}
X_{2.1}=\frac{\partial}{\partial x}.
\end{equation}

  \item[Case 2.2:] $g(t)=kt^b$.\\
  The infinitesimal generators in this case are
  \begin{equation}
X_{2.2.1}=\frac{\partial}{\partial
x},\,\,X_{2.2.2}=(2b-1)x\frac{\partial}{\partial
x}+2t\frac{\partial}{\partial t}+2(b-1)u\frac{\partial}{\partial u}.
\end{equation}

  \item[Case 2.3:] $g(t)=k$.\\
 We obtain the infinitesimal generators as follows
   \begin{equation}
X_{2.3.1}=\frac{\partial}{\partial
x},\,\,X_{2.3.2}=x\frac{\partial}{\partial
x}-2t\frac{\partial}{\partial t}+2u\frac{\partial}{\partial u}.
\end{equation}
\end{description}

\item[Case 3:] $\alpha=\frac{1}{3},\,\ k,b\neq0$.\\
For $\alpha=\frac{1}{3}$, functions of   $g(t)$  can be obtained as
follows
\begin{equation*}
g(t)=k(t-b)^{\frac{2}{3}},\,\,k(t^2-b)^{\frac{1}{3}},\,\,\,ke^{b
t},\,\,kt^b,\, \, k.
\end{equation*}
\begin{description}

\item[Case 3.1:]  $g(t)= k(t-b)^{\frac{2}{3}},\,\,k(t^2-b)^{\frac{1}{3}},\,\,\,ke^{b
t}$.\\
   In these cases, the infinitesimal generator is given by
\begin{equation}
X_{3.1}=\frac{\partial}{\partial x}.
\end{equation}

  \item[Case 3.2:] $g(t)=kt^b$.\\
  The infinitesimal generators in this case are
  \begin{equation}
X_{3.2.1}=\frac{\partial}{\partial
x},\,\,X_{3.2.2}=(3b-1)x\frac{\partial}{\partial
x}+3t\frac{\partial}{\partial t}+(3b-2)u\frac{\partial}{\partial u}.
\end{equation}

  \item[Case 3.3:] $g(t)=k$.\\
 In this case, we obtain the infinitesimal generators as follows
   \begin{equation}
X_{3.3.1}=\frac{\partial}{\partial
x},\,\,X_{3.3.2}=-3t\frac{\partial}{\partial t}+x\frac{\partial}{\partial x}+2u\frac{\partial}{\partial u}.
\end{equation}
\end{description}
\end{description}
\end{The}

\textbf{Proof}. The one-parameter Lie group of transformations in $x,t,u$ with $\varepsilon$ as the group parameter are given
\begin{eqnarray*}
&&t^*=t+\varepsilon\xi^t(t,x,u)+O(\varepsilon^2),\\
&&x^*=x+\varepsilon\xi^x(t,x,u)+O(\varepsilon^2),\\
&&u^*=u+\varepsilon\eta_u(t,x,u)+O(\varepsilon^2),
\end{eqnarray*}
The Lie algebra of K(2,3) equation (Eq. \eqref{ai})  is spanned by vector fields
\begin{eqnarray}\label{aiii}
X = \xi^x(t,x, u) \frac{\partial}{\partial x} + \xi^t(t,x, u)
\frac{\partial}{\partial t} + \eta_u(x, t,
u)\frac{\partial}{\partial u},
\end{eqnarray}
where
\begin{equation}
\xi^x=\frac{dx^*}{d\varepsilon}|_{_{\varepsilon=0}}, \quad
\xi^t=\frac{dt^*}{d\varepsilon}|_{_{\varepsilon=0}}, \quad
\eta_u=\frac{du^*}{d\varepsilon}|_{_{\varepsilon=0}}.
\end{equation}
Applying the $X^{(\alpha)}$ to \eqref{ai}, leads
\begin{equation}\label{aj}
X^{(\alpha)}\left[\frac{\partial^\alpha u}{\partial
t^\alpha}+(u^2)_x+g(t)(u^3)_{xxx}\right]_{\frac{\partial^\alpha u}
{\partial t^\alpha}+(u^2)_x+g(t)(u^3)_{xxx}=0}=0.
\end{equation}

 Expanding the \eqref{aj}, and solving this obtained set using the Maple, we can
distinguish all selections of the function $g(t)$. Finally,
the Lie point symmetries for \eqref{ai} can be obtained as follow.

$\bullet$ If $0<\alpha<1,\,\alpha\neq\frac{1}{2}\,,\frac{1}{3},\
k,b\neq0$, and $g(t)$ be an arbitrary nonvanishing function then we
have:
\begin{eqnarray*}
\xi^x=c_1, \quad \xi^t=0, \quad \eta_u=0.
\end{eqnarray*}
Thus, the infinitesimal generator is given by
\begin{eqnarray*}
X_1=\frac{\partial}{\partial x}.
\end{eqnarray*}

$\bullet$ If $0<\alpha<1,\,\alpha\neq\frac{1}{2}\,,\frac{1}{3},\
k,b\neq0$, and $g(t)=k t^b$ then we have:
\begin{eqnarray*}
\xi^x=c_1+c_2 (\alpha-b) x, \quad \xi^t=- c_2 t, \quad \eta_u=c_2(2
\alpha-b)u.
\end{eqnarray*}
So, the infinitesimal generators are 
\begin{eqnarray*}
X_1=\frac{\partial}{\partial x}, \,\,
X_2=(\alpha-b)x\frac{\partial}{\partial x}-t\frac{\partial}{\partial
t}+(2 \alpha-b)u\frac{\partial}{\partial u}.
\end{eqnarray*}

$\bullet$ If $0<\alpha<1,\,\alpha\neq\frac{1}{2}\,,\frac{1}{3},\
k,b\neq0$, and $g(t)=k$ then we have:
\begin{eqnarray*}
\xi^x=c_1+c_2 \alpha x, \quad \xi^t=- c_2 t, \quad \eta_u=2 c_2
\alpha u.
\end{eqnarray*}
Therefore, the infinitesimal generators are given by
\begin{eqnarray*}
X_1=\frac{\partial}{\partial x}, \,\, X_2=\alpha
x\frac{\partial}{\partial x}-t\frac{\partial}{\partial t}+2 \alpha
u\frac{\partial}{\partial u}.
\end{eqnarray*}
The proof for $\alpha=\frac{1}{2}\,,\frac{1}{3}$ are similar.
Therefore, proof is completed.
%%%%%%%%%%%%%%%%%%%%%%555
\section{Reduced equations and invariant solution of \eqref{ai}} 
  %%%%%%%%%%%%%%%%%%%%%%%%%%%%%%%
	Our purpose for \eqref{ai} is  to reduce it the coordinates $(x,t,u)$ using invariants  $(r,z)$  to a new coordinates\cite{Nadja}.\\
 Let us consider 
\begin{equation*}
X =  \xi^t(t,x, u)
\frac{\partial}{\partial t}+\xi^x(t,x, u)\frac{\partial}{\partial x}  + \eta_u(t,x, u)
\frac{\partial}{\partial u},
\end{equation*}
as a  Lie point symmetry of the time-fractional K(2,3) equation
\begin{equation*}
\frac{\partial^\alpha u}{\partial
t^\alpha}+(u^2)_x+g(t)(u^3)_{xxx}=0, \quad 0<\alpha<1,\,\, t>0.
\end{equation*}
We use two invariants  $z=\psi(x,t)$ and $r=\varphi(x,t)$ which are linearly independent in the characteristic equations
\begin{equation*}
\frac{dt}{\xi^t(t,x,u)}=\frac{dx}{\xi^x(t,x, u)}=\frac{du}{\eta_u(t,x, u)},
\end{equation*}
for  obtaining the invariant solutions. 
After that, we assume one of those invariants is depend to another, 
\begin{equation}\label{300}
z = h(r),
\end{equation}
then we solve \eqref{300} for $u$.
Finally, substituting u  in Eq.\eqref{ai}  for the unknown function $h$, a fractional ODE can be obtained.
Now, we obtain corresponding reduced equations , invariants  and group invariant solutions of \eqref{ai} for different cases of $g(t)$ and $\alpha$ as follows.
\subsection*{{Case1:}}
\begin{description}
  \item[$\bullet$ Case 1.1:]$0<\alpha<1$,\,$\alpha\neq\frac{1}{2}\,,\,\frac{1}{3}$ and $g(t)$ is a  nonvanishing  arbitrary function.
  \item[$\bullet$ Case 1.2:]\,$\alpha=\frac{1}{2}$, $g(t)=\{k, \,\,kt^b, \,\,
  ke^{bt}\}$
  \item[$\bullet$ Case 1.3:]\,$\alpha=\frac{1}{3}$, $g(t)=\{k, \,kt^b, \,
  ke^{bt},\,k(t-b)^{\frac{2}{3}},\,k(t^2-b)^{\frac{1}{3}}\}$
\end{description}
In these cases, according to the infinitesimal generator
$X=\frac{\partial}{\partial x}$, the similarity
  variables using the method of characteristics   are as follows:
\begin{equation}
 z= u, \,\,r=t,
\end{equation}
and a solution is
\begin{equation}\label{invf}
z=h(r)\Rightarrow  u=h(t).
\end{equation}
By substituting \eqref{invf}  into \eqref{ai} we find the $h(r)$.
Thus $h(r)$ must be satisfied:
\begin{equation}\label{akb}
\frac{d^\alpha h(t)}{dt^\alpha}=0.
\end{equation}
Then by solving the above equation by the Laplace transform\cite{Pod}, we have
\begin{equation}
h(t)=\frac{\kappa t^{\alpha-1}}{\Gamma(\alpha)},\quad \kappa \,\text{ is a constant}.
\end{equation}

\subsection*{{Case2:}}
\begin{description}
  \item[$\bullet$ Case 2.1:]$\alpha\neq\frac{1}{2}, \,\frac{1}{3},\, g(t)=kt^b$.\\
In this case
  \begin{equation}
X_{1.2.2}=-t\frac{\partial}{\partial t}+(\alpha-b)x\frac{\partial}{\partial
x}+(2\alpha-b)u\frac{\partial}{\partial
u},
\end{equation}
so  the similarity   variables for this Lie point symmetry  using
the method of characteristics are as follows:
\begin{equation}
  r=tx^{\frac{1}{\alpha-b}}, \quad z=u x^{\frac{b-2\alpha}{\alpha-b}},
\end{equation}
and a solution  for equation \eqref{ai} is
\begin{equation}\label{invys}
z=h(r)\Rightarrow  u= x^{\frac{b-2\alpha}{b-\alpha}} h(tx^{\frac{1}{\alpha-b}}).
\end{equation}
We substitute \eqref{invys}  into \eqref{ai} to find the
$f(r)$ and $f(r)$ must be satisfied in the fractional ODE as follows:

\begin{eqnarray*}
&&(b-\alpha)^3\frac{\partial^\alpha h}{\partial
r^\alpha}-2(b-\alpha)^2r h(r) h'(r)+18(-1+2b-5\alpha)kr^{2+b} h(r)
h'(r)^2\\
&&-18kr^{3+b}h(r)h'(r)h''(r)+3(2b^3-17b^2\alpha+46b\alpha^2-40\alpha^3)kr^bh(r)^3\\
&&-6kr^{3+b}h'(r)^3+(2b^3-8b^2\alpha+10b\alpha^2-4\alpha^3)h(r)^2-3(11b^2+15\alpha+74\alpha^2\\
&&-2b(29\alpha+3)+1)kr^{1+b}h'(r)h(r)^2
+9(-1+2b-5\alpha)kr^{2+b}h''(r)h(r)^2\\
&&-3kr^{3+b}h'''(r)h(r)^2=0.
\end{eqnarray*}
Where $\alpha\neq\frac{1}{2}, \,\frac{1}{3}$.
\end{description}

\begin{description}
  \item[$\bullet$ Case 2.2:]$\alpha\neq\frac{1}{2}, \,\frac{1}{3},\, g(t)=k$.\\
For this case we have
  \begin{equation}
X_{1.3.2}=\alpha x\frac{\partial}{\partial
x}-t\frac{\partial}{\partial t}+2\alpha u\frac{\partial}{\partial
u},
\end{equation}
so  the similarity   variables for this Lie point symmetry  using
the method of characteristics are as follows:
\begin{equation}
  r=tx^{\frac{1}{\alpha}}, \quad z=u x^{-2},
\end{equation}
and a solution for \eqref{ai}  is
\begin{equation}\label{invy5}
z=h(r)\Rightarrow  u= x^2 h(tx^{\frac{1}{\alpha}}).
\end{equation}
Again we substitute \eqref{invy5}  into \eqref{ai} to obtain the
$f(r)$. So $h(r)$ must satisfy in the fractional ODE as follows:
\begin{eqnarray*}
&&\alpha^3\frac{\partial^\alpha h}{\partial r^\alpha}+2\alpha^2 r
h(r) h'(r)
+18(5\alpha+1)kr^{2}h(r) h'(r)^2\\
&&+18kr^{3}h(r) h'(r) h''(r)+120k\alpha^3h(r)^3+6kr^3
h'(r)^3+4\alpha^3 h(r)^2\\
&&+3(74\alpha^2+15\alpha+1)k r h(r)^2 h'(r)+9(5\alpha+1)kr^{2}h(r)^2 h''(r)\\
&&+3kr^{3} h(r)^2 h'''(r)=0 .
\end{eqnarray*}
Where $\alpha\neq\frac{1}{2}, \,\frac{1}{3}$.
\end{description}

\subsection*{{Case3:}}

\begin{description}
  \item[$\bullet$ Case 3.1:]$\alpha=\frac{1}{2}, \, g(t)=kt^b$.\\
For this case we have
  \begin{equation}
X_{2.2.2}=(2b-1)x\frac{\partial}{\partial
x}+2t\frac{\partial}{\partial t}+2(b-1)u\frac{\partial}{\partial u},
\end{equation}
so  the similarity   variables for this Lie point symmetry  using
the method of characteristics are as follows:
\begin{equation}
  r=tx^{\frac{2}{1-2   b}}, \quad z=u x^{\frac{2 b-2}{1-2 b}},
\end{equation}
and in view of \eqref{invf}, a solution for \eqref{ai}  is
\begin{equation}\label{invy}
u= x^{\frac{2-2 b}{1-2 b}}
h\left(tx^{\frac{2}{1-2   b}}\right).
\end{equation}
We substitute \eqref{invy}  into \eqref{ai} to obtain $h(r)$.
After that  $h(r)$ must be satisfied in the FDE as follows:

\begin{eqnarray*}
&&\frac{(2b-1)^3}{4}\frac{\partial^\alpha h}{\partial
r^\alpha}-(1-2b)^2r h(r) h'(r)+18(4b-7)kr^{2+b} h(r)
h'(r)^2\\
&&-36kr^{3+b}h(r)
h'(r)h''(r)+3(4b^3-17b^2+23b-10)kr^bh(r)^3\\
&&-12kr^{3+b}h'(r)^3+(4b^3-8b^2+5b-1)h(r)^2-6(11b^2-35b+27)kr^{1+b}h'(r)h(r)^2\\
&&+9(4b-7)kr^{2+b}h''(r)h(r)^2-6kr^{3+b}h'''(r)h(r)^2=0.
\end{eqnarray*}
Where $\alpha=\frac{1}{2}$.
\end{description}

\begin{description}
  \item[$\bullet$ Case 3.2:]$\alpha=\frac{1}{2}, \, g(t)=k$.\\
For this case we have
  \begin{equation}
X_{2.3.2}=x\frac{\partial}{\partial x}-2t\frac{\partial}{\partial
t}+2u\frac{\partial}{\partial u},
\end{equation}
so  the similarity   variables for this Lie point symmetry  using
the method of characteristics are as follows:
\begin{equation}
  r=tx^2, \quad z=u x^{-2},
\end{equation}
and in view of \eqref{invf}, a solution for \eqref{ai}  is
\begin{equation}\label{invy}
u= x^2 h\left(tx^2\right).
\end{equation}
To obtain $f(r)$, we substitute \eqref{invy}  into \eqref{ai}.
Then $f(r)$ must satisfy in the fractional ODE as follows:

\begin{eqnarray*}
&&\frac{1}{4}\frac{\partial^\alpha f}{\partial r^\alpha}+30kf(r)^3
+12kr^{3}f'(r)^3+r f(r) f'(r)+126kr^{2}f(r) f'(r)^2\\
&&+36kr^{3}f(r) f'(r) f''(r)+f(r)^2+162kr
f'(r)f(r)^2+63kr^{2}f(r)^2 f''(r)\\
&&+6kr^{3} f(r)^2 f'''(r)=0 .
\end{eqnarray*}
Where $\alpha=\frac{1}{2}$.
\end{description}

\subsection*{{Case4:}}
\begin{description}
  \item[$\bullet$ Case 4.1:]$\alpha=\frac{1}{3}, \, g(t)=kt^b$.\\
For this case we have
  \begin{equation}
X_{3.2.2}=(3b-1)x\frac{\partial}{\partial
x}+3t\frac{\partial}{\partial t}+(3b-2)u\frac{\partial}{\partial u},
\end{equation}
so  the similarity   variables for this Lie point symmetry  using
the method of characteristics are as follows:
\begin{equation}
  r=tx^{\frac{3}{1-3b}}, \quad z=u x^{\frac{3 b-2}{1-3 b}},
\end{equation}
and a solution to our equation  is
\begin{equation}\label{invy0}
z=h(r)\Rightarrow  u= x^{\frac{2-3 b}{1-3 b}}
h\left(tx^{\frac{3}{1-3 b}}\right).
\end{equation}
We substitute \eqref{invy0}  into \eqref{ai} to determine the
$h(r)$. Then $h(r)$ must satisfy in the fractional ODE as follows:

\begin{eqnarray*}
&&{(3b-1)^3}\frac{\partial^\alpha h}{\partial r^\alpha}-6(1-3b)^2r
h(r) h'(r)+324(3b-4)kr^{2+b} h(r)
h'(r)^2\\
&&-486kr^{3+b}h(r)
h'(r)h''(r)+3(54b^3-153b^2+138b-40)kr^bh(r)^3\\
&&-162kr^{3+b}h'(r)^3+(54b^3-72b^2+30b-4)h(r)^2-9(99b^2-228b+128)kr^{1+b}h'(r)h(r)^2\\
&&+162(3b-4)kr^{2+b}h''(r)h(r)^2-81kr^{3+b}h'''(r)h(r)^2=0.
\end{eqnarray*}
Where $\alpha=\frac{1}{3}$.
\end{description}

\begin{description}
  \item[$\bullet$ Case 4.2:]$\alpha=\frac{1}{3}, \, g(t)=k$.\\
For this case we have
  \begin{equation}
X_{3.3.2}=x\frac{\partial}{\partial x}-3t\frac{\partial}{\partial
t}+2u\frac{\partial}{\partial u},
\end{equation}
so  the similarity   variables for this Lie point symmetry  using
the method of characteristics are as follows:
\begin{equation}
  r=tx^3, \quad z=u x^{-2},
\end{equation}
and a solution to our equation  is
\begin{equation}\label{invym}
z=h(r)\Rightarrow  u= x^2 h\left(tx^3\right).
\end{equation}
We substitute \eqref{invym}  into \eqref{ai} to determine the
$h(r)$. Then $h(r)$ must satisfy in the fractional ODE as follows:

\begin{eqnarray*}
&&\frac{\partial^\alpha h}{\partial r^\alpha}+120kh(r)^3
+162kr^{3}h'(r)^3+6r h(r) h'(r)+1296kr^{2}h(r) h'(r)^2\\
&&+486kr^{3}h(r) h'(r) h''(r)+4h(r)^2+1152kr
h'(r)h(r)^2+648kr^{2}h(r)^2h''(r)\\
&&+81kr^{3} h(r)^2 h'''(r)=0 .
\end{eqnarray*}
Where $\alpha=\frac{1}{3}$.
\end{description}

\section{Conclusion}

The corresponding invariants and group invariant solutions has been obtained for  the time-fractional K(m,n) equation when the the fractional derivative is in the Riemann-Liouville sense. Finally,  we reduced this time-fractional  equation into a nonlinear ODE of fractional order. For this propose, the Lie group method  and the symmetry properties have been investigated for the governing equation. They have been used to reduced the give FPDE to a corresponding FDE which might be solved easily.

\section*{Conflict of interest}
 The author declare that they have no conflict of interest.
\section*{Funding}
Not applicable.
\section*{Authors' contributions}
All authors have read and approved the final manuscript.
\section*{Acknowledgments}
I would like thank the editor and reviewers to consider this paper for review.
%\begin{acknowledgements}

%\section*{Acknowledgements}
%\noindent  The authors would like to extend great gratitude to the Editor and anonymous reviewers, whose insightful comments and constructive suggestions helped us to significantly improve the quality of this paper.
%\newpage

\end{document}